\documentclass[reqno,10pt]{amsart}    

\usepackage{amscd,amsfonts,mathrsfs,amsthm,enumerate}
\usepackage{amssymb, amsmath}

\usepackage{color}           
\definecolor{myowncolor1}{rgb}{0,0,.75}
\definecolor{darkred}{rgb}{0.7,0,0}
\definecolor{darkgreen}{rgb}{0,0.5,0}
\definecolor{darkblue}{rgb}{0,0,.75}
\definecolor{orange}{rgb}{1,.5,0}
\definecolor{papier}{rgb}{1,1,.97}

\def\complex  #1{{\mathbb C^{#1}}}

\def\cont     {\lrcorner}

\def\nb       {\overline{\nabla}}

\def\picture#1#2{
\begin{latexonly}
\ifx\pdfoutput\undefined 
    \includegraphics[width=#2\hsize]
      {#1.eps}
\else
    \includegraphics[width=#2\hsize]
      {#1.jpg}
\fi
\end{latexonly}
}

\newtheorem{theorem}{Theorem}[section]   
\newtheorem{lemma}[theorem]{Lemma}   
\newtheorem{corollary}[theorem]{Corollary}   
\newtheorem{proposition}[theorem]{Proposition}   
\newtheorem{remark}[theorem]{Remark}   
  
\theoremstyle{definition}   
   
\newtheorem{example}[theorem]{Example}   
   
\numberwithin{equation}{section}

\newcommand{\bfig}{\begin{figure}}
\newcommand{\efig}{\end{figure}}

{}

\begin{document}
\title
[Lagrangian submanifolds in nearly K\"ahler manifolds]
{Decomposition and minimality of Lagrangian submanifolds in nearly K\"ahler manifolds}

\author[Sch\"afer]{\sc Lars Sch\"afer\,}
\email{schaefer@math.uni-hannover.de} 
\author[Smoczyk]{\sc \,Knut Smoczyk}
\email{smoczyk@math.uni-hannover.de}
\address{Leibniz Universit\"at Hannover, Institut f\"ur Differentialgeometrie,    
  Welfengarten 1, 30167 Hannover, Germany}

\begin{abstract}
We show that 
Lagrangian submanifolds in 
six-dimensional nearly K\"ahler (non K\"ahler) manifolds and in twistor spaces
$Z^{4n+2}$ over quaternionic K\"ahler manifolds $Q^{4n}$ are minimal. Moreover, we
will prove that any Lagrangian submanifold $L$ in a nearly K\"ahler manifold
$M$ splits into a product of two Lagrangian submanifolds for which
one factor is Lagrangian in the strict nearly K\"ahler part of $M$ and the second
factor is Lagrangian in the K\"ahler part of $M$. Using this splitting theorem
we then describe Lagrangian submanifolds in nearly K\"ahler manifolds
of dimensions six, eight and ten.
\end{abstract}

\date{\today}   
\thanks{Supported by DFG, priority program SPP 1154, SM 78/4-2}
\maketitle

\section{Introduction}
Suppose $(M,\omega)$ is a symplectic manifold. Then a submanifold $L\subset M$ is called
{\it Lagrangian}, if $\omega_{|TL}=0$ and $2\dim(L)=\dim M$. If $L$ can be embedded
(immersed) as a Lagrangian submanifold in $\complex{n}$, then by Darboux's theorem $L$ can 
also be embedded (immersed) into any other symplectic manifold $(M,\omega)$. One
of the most interesting examples for symplectic manifolds are K\"ahler manifolds
$(M,J,g)$, where the symplectic form $\omega$ is given by $\omega(X,Y)=g(JX,Y)$. 

In this paper we will consider Lagrangian submanifolds in nearly K\"ahler manifolds.
An almost Hermitian manifold $(M,J,g)$ is called 
{\it nearly K\"ahler}, if its almost complex structure $J$ satisfies
\begin{equation}\label{nk1}
\nabla_X(J)X=0\,,\quad\forall\,X\in TM\,.
\end{equation}
A nearly K\"ahler manifold $(M,J,g)$ is K\"ahler, iff $\nabla J=0$. $(M,J,g)$ is called
strict nearly K\"ahler, if $\nabla_XJ\neq 0$ for all $X\in TM$, $X\neq 0$.

In contrast to the K\"ahler case, nearly K\"ahler (non K\"ahler) manifolds are not
symplectic manifolds (at least not with
their characteristic $2$-form $\omega$).

Anyway, also in this case, it is common 
to say that $L\subset M$ is Lagrangian, if $\omega_{|TL}=0$ and if $L$ has half the
dimension of $M$. Under the weaker condition that $(M,J,g)$ is nearly K\"ahler,
it may not be possible to 
find Lagrangian submanifolds of a given topological type $L$ in $M$, 
even if $L$ is a Lagrangian submanifold in $\complex{n}$. This is because there
is no analogue of Darboux's theorem in the nearly K\"ahler case. Nevertheless, there 
exist numerous examples for nontrivial Lagrangian submanifolds in the nearly
K\"ahler six-sphere $S^6$ \cite{Vrancken-2003}.

Nearly K\"ahler geometry was first studied  in the 1970s by Gray 
(see \cite{Gray-1970},\cite{Gray-1976}) in the context of weak holonomy. 
The most prominent example of a nearly K\"ahler manifold is $S^6$ with its
standard almost complex structure and Riemannian metric. We shortly resume the 
state of the classification of nearly K\"ahler manifolds $M$. The reader is
invited to consult section 5 of \cite{Butruille-2006} for a short survey with detailed
information and more references. In the rest of this paragraph we consider $M$
to be complete and simply connected. After splitting of the K\"ahler factor
$M$ can supposed to be strict nearly K\"ahler. Nagy  \cite{Nagy-2002a}
has reduced the classification of strict nearly K\"ahler manifolds using
previous work of Cleyton and Swann \cite{Cleyton-Swann} to 
almost hermitian products of:
\begin{itemize}
\item[(i)] (naturally reductive) three-symmetric spaces \cite{Gray-1972, Gray-Wolf}, 
\item[(ii)] twistor spaces of non locally symmetric positive quaternionic K\"ahler
manifolds  endowed with their non-integrable complex structure,
\item[(iii)] irreducible six-dimensional nearly K\"ahler manifolds. 
\end{itemize}
The class of examples from three-symmetric spaces can be divided into three types.
One type are nearly K\"ahler structures on twistor spaces over symmetric
spaces. In this case the isotropy representation is complex reducible. Another
type are those with irreducible isotropy representation and for the third type
the  isotropy representation is real reducible.
In the case of a homogenous  nearly K\"ahler six-manifold
 Butruille \cite{Butruille-2005} showed that there 
are only the previous known examples \cite{Gray-1972, Gray-Wolf}
\begin{eqnarray*}
S^6&=&G_2/SU_3, S^3 \times S^3=SU_2\times SU_2\times SU_2/SU_2,\\ \mathbb
 CP^3&=&SO_5/U_2\times S^1 \mbox{ and } F(1,2)=SU_3/T^2,
\end{eqnarray*}
 which are all three-symmetric spaces with the metric defined by the Killing
 form.
 In dimension eight strict nearly K\"ahler manifolds are almost hermitian products of a
six dimensional nearly K\"ahler factor and a Riemannian surface. In dimension
ten  strict nearly K\"ahler manifolds are either almost hermitian products of a
six dimensional nearly K\"ahler factor and a complex surface or they come from
twistor spaces. Therefore we later study twistor spaces and three-symmetric  spaces in more detail.

Since there is no Darboux theorem, the existence of Lagrangian submanifolds in nearly
K\"ahler manifolds, even locally, is not unobstructed. Indeed, in this paper we will show 
that Lagrangian submanifolds tend to be minimal in the strict nearly K\"ahler directions
of $M$, a fact which was previously known only for the nearly K\"ahler six-sphere $S^6$ \cite{Ejiri-1981}. 
To be precise, among other things we will prove the following three theorems:

{\bf Theorem A:}
{\it
If $L$ is a Lagrangian submanifold in a strict nearly K\"ahler six-manifold, 
then $L$ is orientable and minimal.}

{\bf Theorem B:}
{\it
If $L\subset Z^{4n+2}$ is a Lagrangian submanifold in a twistor space over a positive quaternionic
K\"ahler manifold $M^{4n}$, then $L$ is minimal. Moreover if $n>1$, then 
$TL=\mathcal{D}^\perp\oplus\mathcal{D}$ decomposes into a $2n$-dimensional horizontal
distribution $\mathcal{D}^\perp$ and a vertical line bundle $\mathcal{D}$ and the
second fundamental form $II(\mathcal{D},\cdot)$ of the vertical part vanishes identically.
}

{\bf Theorem C:}
{\it Let $M$ be a nearly K\"ahler manifold and $L\subset M$ be Lagrangian. Then $M$ and $L$
decompose locally into products $M=M_K\times M_{SNK}$, $L=L_K\times L_{SNK}$, where
$M_K$ is K\"ahler, $M_{SNK}$ is strict nearly K\"ahler and $L_K\subset M_K$, $L_{SNK}\subset M_{SNK}$
are both Lagrangian. The dimension of $L_K$ is given by
$$\dim L_K=\frac{1}{2}\dim\operatorname{ker}(r)\,,$$
where $r=\operatorname{Ric}-\operatorname{Ric}^*$ is defined as in (\ref{op r}).
Moreover, if the splitting of $M$ is global and $L$ is simply connected,
then $L$ decomposes globally as well.
}

The organization of the paper is as follows: 
In section \ref{sec 1} we will recall the basic geometric quantities related to
nearly K\"ahler geometry. Section \ref{sec 2} concerns the geometry of Lagrangian submanifolds
$L$ in nearly K\"ahler manifolds. Here we will define some new tensors, that are couplings
of the second fundamental form $II$ of $L$ with the torsion of the
canonical connection $\nb$. Then we will prove some important identities related to
these quantities. In section \ref{sec 3} we will consider the case of Lagrangian
submanifolds in six-dimensional nearly K\"ahler manifolds and we will prove our
main theorem A. Section \ref{sec 5} is devoted to the decomposition of Lagrangian
submanifolds $L$ in nearly K\"ahler manifolds $M$ and we will show that the Lagrangian 
condition is compatible with the splitting of the nearly K\"ahler manifold $M$ induced by the
operator $r=\operatorname{Ric}-\operatorname{Ric}^*$. This will be theorem C. We will
need this result to prove theorem B in the following section, where we analyse Lagrangian
submanifolds in twistor spaces over positive quaternionic K\"ahler manifolds in more
detail. In the remainder we will prove that left invariant Lagrangian submanifolds
in three-symmetric spaces are totally geodesic (Section \ref{sec 7}) and then we
will analyse in Section \ref{sec 8} the space of deformations of a Lagrangian submanifold in strict 
nearly K\"ahler six-manifolds in terms of coclosed eigenforms of the Hodge-Laplacian.

\section{The geometry of nearly K\"ahler manifolds}\label{sec 1}

Let $(M,J,g)$ be a smooth nearly K\"ahler manifold with Riemannian metric
$\langle\cdot,\cdot\rangle:=g(\cdot,\cdot)$, Levi-Civita connection $\nabla$
and almost complex structure $J\in\operatorname{End}(TM)=\Gamma(TM\otimes T^*M)$.
Then the nearly K\"ahler condition is
\begin{eqnarray}
\langle JX,JY\rangle&=&\langle X,Y\rangle\,,\quad\forall\, X,\,Y\in TM\,,\nonumber\\
J^2&=&-\operatorname{Id}\,,\nonumber\\
\nabla_X(J)Y&=&-\nabla_Y(J)X\,,\quad\forall\, X,\,Y\in TM\,.\nonumber
\end{eqnarray}
This gives
\begin{equation}\label{nk4}
\langle\nabla_X(J)Y,Z\rangle=-\langle\nabla_Y(J)X,Z\rangle=-\langle\nabla_X(J)Z,Y\rangle
\end{equation}
and
\begin{equation}\label{domega}
d\omega=3\nabla\omega\,.
\end{equation}
Let 
$$R(W,X,Y,Z)=\langle\nabla_{[W,X]}Y,Z\rangle-\langle[\nabla_W,\nabla_X]Y,Z\rangle$$
denote the value of the curvature tensor of $M$ on vector fields $W,X,Y,Z\in\mathscr{X}(M)$.
The following is well known (see \cite{Gray-1970},\cite{Gray-1976}):
\begin{eqnarray}
\nabla_{JX}(J)Y&=&\nabla_{X}(J)JY=-J\nabla_X(J)Y\,.\label{nk5}
\end{eqnarray}
Let $\{e_1,\dots,e_{2n}\}$ be a local orthonormal frame field. Then the Ricci and
Ricci$\phantom{}^*$ curvatures are
$$\operatorname{Ric}(X,Y)=\sum_{i=1}^{2n}R(X,e_i,Y,e_i)\,,$$
$$\operatorname{Ric}^*(X,Y)=\frac{1}{2}\sum_{i=1}^{2n}R(X,JY,e_i,Je_i)\,.$$
We define the endomorphism
\begin{equation}\label{op r}
r\in\operatorname{End}(TM)\,,\quad\langle rX,Y\rangle:=\operatorname{Ric}(X,Y)-\operatorname{Ric}^*(X,Y)\,.
\end{equation}
Then $r$ satisfies the following identities (see \cite{Koto-1960}):
\begin{eqnarray}
[r,J]&=&0\,,\label{koto 1}\\
\langle rX,Y\rangle&=&\sum_{i=1}^{2n}\langle\nabla_X(J)e_i,\nabla_Y(J)e_i\rangle
=\sum_{i=1}^{2n}\langle\nabla_{e_i}(J)X,\nabla_{e_i}(J)Y\rangle\,,\label{koto 2}\\
r&=&-\sum_{i=1}^{2n}\nabla^2_{e_ie_i}(J)J\,.\label{koto 3}
\end{eqnarray}
On a nearly K\"ahler manifold we can define another connection $\nb$ by
$$\nb_XY:=\nabla_XY-\frac{1}{2}J\nabla_X(J)Y\,.$$
As has been shown in \cite{Friedrich-Ivanov}, $\nb$ is the uniquely defined connection
on $M$ with
\begin{eqnarray}
\nb g&=&0\,,\label{nb g}\\
\nb J&=&0\,\label{nb J}
\end{eqnarray}
such that 
\begin{equation}\label{tau}
\tau(X,Y,Z):=g(T(X,Y),Z)
\end{equation}
is completely skew symmetric, where $T(X,Y)$ denotes the torsion of $\nb$. 
This can be compared to the
Bismut connection in complex geometry. In our case
the torsion $T(X,Y)$ now becomes
$$T(X,Y)=-J\nabla_X(J)Y\,.$$
By the nearly K\"ahler condition, the torsion satisfies
\begin{equation}\label{tau2}
T(JX,Y)=-JT(X,Y)=T(X,JY)\,.
\end{equation}

It has been shown in \cite{Belgun-Moroianu} that $\nabla(J)$ and the torsion $T$
are $\nb$-parallel. 
Nagy proved the following theorem in \cite{Nagy-2002}:

{\bf Theorem: }(Nagy)
{\it
Any complete, strict nearly K\"ahler manifold $M$ has positive Ricci curvature and
hence is compact with finite fundamental group.}

The operator $r$ is also parallel w.r.t. $\nb$, i.e.
$\nb r=0$ (cf. \cite{Nagy-2002}). 
In fact in this reference it is
further proven that even the  Ricci and Ricci$\phantom{}^*$ curvatures are
$\nb$-parallel. Since $r$ is parallel w.r.t. a metric connection $\nb$, the spectrum of $r$ does not depend
on $p\in M$. Moreover, $r$ is selfadjoint and positive semidefinite. Note, that
$(M,J,g)$ is strict nearly K\"ahler, iff $r$ is positive definite, which is
equivalent to $\operatorname{ker}(r)=\{0\}$. In general, since $r$ commutes with $J$,
all eigenspaces $E_p(\lambda)\subset T_pM$ of $r$ are complex subspaces of $T_pM$
and in particular they are even dimensional.

\section{Lagrangian submanifolds}\label{sec 2}

{F}or the rest of this section let us assume that $L\subset M$ is a Lagrangian 
submanifold of a nearly K\"ahler manifold $(M^{2n},J,g)$. Then since $n=\dim(L)=\frac{1}{2}\dim(M)$
we have
$$\langle JX,Y\rangle=0\,,\quad\forall\,X,\,Y\in TL\quad\Leftrightarrow\quad
J:TL\to T^\perp L\quad\text{is an isomorphism}.$$
{F}rom $\omega_{|TL}=0$ we deduce $d\omega_{|TL}=0$. On the other hand (\ref{domega})
implies
$$d\omega(X,Y,Z)=3\langle\nabla_X(J)Y,Z\rangle\,.$$
{F}rom this and (\ref{nk4}) the following Lemma easily follows (see also \cite{Hou}).
\begin{lemma}\label{lemma 1}
Suppose $L\subset M$ is a Lagrangian submanifold in a nearly K\"ahler manifold $(M,J,g)$.
Then
\begin{eqnarray}
&&\nabla_X(J)Y\in T^\perp L\,,\quad\forall\,X,\,Y\in TL\,,\label{lag 1}\\
&&\nabla_X(J)Y\in T^\perp L\,,\quad\forall\,X,\,Y\in T^\perp L\,,\label{lag 2}\\
&&\nabla_X(J)Y\in TL\,,\text{ if }\,X\in TL, Y\in T^\perp L\,\text{ or if }\,X\in T^\perp L, Y\in TL\,.\label{lag 3}
\end{eqnarray}
\end{lemma}
Denote by $II$ the second fundamental form of the Lagrangian immersion $ L \subset M^{2n}$
into a nearly K\"ahler manifold $M.$ 
\begin{proposition} \label{2nd_fund_info}
{F}or a Lagrangian submanifold $ L \subset M^{2n}$ in a nearly K\"ahler manifold we have the
following information.
\begin{itemize}
\item[(i)]
The second fundamental form is given by $\langle II(X,Y),U \rangle = \langle
\nb_XY,U \rangle$ for $X,Y \in TL$ and $U \in
T^\perp L.$
\item[(ii)] The tensor 
$C(X,Y,Z):=\langle II(X,Y),JZ \rangle =\omega(II(X,Y),Z)$, $\forall\,X,Y,Z \in TL$ is totally symmetric.
\end{itemize}
\end{proposition}
\begin{proof}
{F}rom Lemma \ref{lemma 1} we compute for  $X,Y \in TL$ and $U \in
T^\perp L$ the second fundamental form  $II$ 
$$\langle II(X,Y),U \rangle  = \langle \nabla_XY,U \rangle = \langle  \nb_X Y -\frac{1}{2} J\nabla_X (J)Y,U \rangle = \langle \nb_XY,U \rangle\,.$$ 
This yields part (i). Next we prove (ii): First we
  observe for $X,Y,Z \in TL$ 
\begin{eqnarray*} 
C(X,Y,Z)&=&\langle II(X,Y),JZ \rangle = \langle \nb_XY,JZ \rangle = - \langle Y,\nb_X(JZ) \rangle\\ &=& -\langle Y,J\nb_X
Z \rangle=\langle \nb_XZ,JY \rangle =C(X,Z,Y).
\end{eqnarray*}
Since the second fundamental form is symmetric, it follows that $C$ is totally symmetric.
\end{proof}

Next we generalize an identity of \cite{Ejiri-1981} to nearly K\"ahler
manifolds of arbitrary dimension. This and the next lemma will be crucial to prove
that Lagrangian submanifolds in strict nearly K\"ahler six-manifolds and in twistor
spaces $Z^{4n+2}$ over quaternionic K\"ahler manifolds with their canonical nearly
K\"ahler structure are minimal.
\begin{lemma} \label{Lemma cyclic id}
The second fundamental form $II$ of a Lagrangian immersion $L \subset M^{2n}$
into a nearly K\"ahler manifold and the tensor $\nabla(J)$ satisfy the following identity
\begin{equation}
\quad \quad \langle II( X,J\nabla_Y(J)Z),U \rangle = \langle J\nabla_{
  II(X,Y)}(J)Z,U \rangle + \langle J\nabla_Y(J) II(X,Z),U \rangle
\end{equation}
with $X,Y \in TL$ and $U \in T^\perp L.$
\end{lemma}
\begin{proof}
The proof of this identity uses $\nb (J)=0,$ $\nb (\nabla
J)=0$ and Lemma \ref{lemma 1}. With $X,Y,Z \in TL$ and $U \in  T^\perp L$ we
obtain 
\begin{eqnarray*}
\langle II(X,J\nabla_Y(J)Z),U \rangle &=& \langle \nb_X(J\nabla_Y(J)Z),U
\rangle \overset{\nb (J)=0}=
\langle J\nb_X(\nabla_Y(J)Z),U \rangle\\
&\overset{\nb (\nabla (J))=0}=&\langle J \left[ \nabla_{\nb_XY}(J)Z+\nabla_Y(J)\nb_X Z \right],U \rangle \\
&=&\langle J[ -\nabla_Z(J)\nb_X Y + \nabla_Y(J) \nb_X Z ],U \rangle\\
&=& -\langle {\nb _XY}, {\nabla_Z(J) J
  U}\rangle +  \langle {\nb_XZ},{\nabla_Y(J)JU}\rangle\\
&=& -\langle {II(X,Y)}, {\nabla_Z(J) J
  U}\rangle +  \langle {II(X,Z)},{\nabla_Y(J)JU}\rangle\\
&=& -\langle J\nabla_Z(J) II(X,Y),U\rangle+ \langle J\nabla_Y(J) II(X,Z),U\rangle \\
&=& \langle J\nabla_{ II(X,Y)}(J)Z,U\rangle + \langle J\nabla_Y(J) II(X,Z),U\rangle.
\end{eqnarray*}
This is exactly the claim of the Lemma. 
\end{proof}
Given the tensor $C(X,Y,T(Z,V))$ we define the following traces
\begin{eqnarray}
\alpha(X,Y)&:=&\sum_{i=1}^nC(e_i,X,T(e_i,Y))\,,\nonumber\\
\beta(X,Y,Z)&:=&\sum_{i=1}^nC(T(e_i,X),Y,T(e_i,Z))\,,\nonumber
\end{eqnarray}
where $\{e_1,\dots,e_n\}$ is a local orthonormal frame of $TL$. 
\begin{lemma}\label{lemma cyc}
{F}or a Lagrangian immersion in a nearly K\"ahler manifold and any $X,Y,Z,V\in TL$ holds:
\begin{eqnarray}
&&C(X,Y,T(Z,V))+C(X,Z,T(V,Y))+C(X,V,T(Y,Z))=0\,,\label{cyc 2}\\
&&\alpha(X,Y)-\alpha(Y,X)=\langle\overrightarrow H,JT(X,Y)\rangle\,,\label{cyc mean}\\
&&\beta(X,Y,Z)=\beta(Z,Y,X)=\beta(Y,X,Z)+\alpha(T(Y,X),Z)\,,\label{cyc beta}\\
&&\alpha(T(X,Y),Z)+\alpha(T(Y,Z),X)+\alpha(T(Z,X),Y)=0\,.\label{cyc beta2}
\end{eqnarray}
Here $\overrightarrow H$ denotes the mean curvature vector of $L$.
\end{lemma}
\begin{proof}
Let us first rewrite the identity in Lemma \ref{Lemma cyclic id} in terms of the tensor $C$
and the torsion $T(X,Y)=-J\nabla_X(J)Y$. Let $X,Y,Z,V\in TL$ be arbitrary. Then Lemma
\ref{Lemma cyclic id} gives
\begin{eqnarray}
C(X,T(Z,Y),V)&=&\langle II( X,J\nabla_Y(J)Z),JV \rangle\nonumber\\ 
&=& \langle J\nabla_{II(X,Y)}(J)Z,JV \rangle + \langle J\nabla_Y(J) II(X,Z),JV \rangle\nonumber\\
&=&\langle J\nabla_Z(J)(JV),II(X,Y)\rangle-\langle J\nabla_Y(J)(JV),II(X,Z)\rangle\nonumber\\
&=&-\langle J^2\nabla_Z(J)V,II(X,Y)\rangle+\langle J^2\nabla_Y(J)V,II(X,Z)\rangle\nonumber\\
&=&\langle JT(Z,V),II(X,Y)\rangle-\langle JT(Y,V),II(X,Z)\rangle\nonumber\\
&=&C(X,Y,T(Z,V))-C(X,Z,T(Y,V))\,.\nonumber
\end{eqnarray}
This is (\ref{cyc 2}). Taking a trace gives
\begin{eqnarray}
\alpha(X,Y)&=&\sum_{i=1}^nC(e_i,X,T(e_i,Y))\nonumber\\
&\overset{(\ref{cyc 2})}{=}&\sum_{i=1}^nC(e_i,e_i,T(X,Y))+\sum_{i=1}^nC(e_i,Y,T(e_i,X))\nonumber\\
&=&\langle\overrightarrow H,JT(X,Y)\rangle+\alpha(Y,X)\,,\nonumber
\end{eqnarray}
which is (\ref{cyc mean}). The first identity in (\ref{cyc beta}) is clear since $C$ is
fully symmetric. If we apply (\ref{cyc 2}) to $\beta(X,Y,Z)$, then we get
\begin{eqnarray}
\beta(X,Y,Z)
&=&\sum_{i=1}^nC(T(e_i,X),Y,T(e_i,Z))\nonumber\\
&=&-\sum_{i=1}^nC(T(X,Y),e_i,T(e_i,Z))-\sum_{i=1}^nC(T(Y,e_i),X,T(e_i,Z))\nonumber\\
&=&\sum_{i=1}^nC(e_i,T(Y,X),T(e_i,Z))+\sum_{i=1}^nC(T(e_i,Y),X,T(e_i,Z))\nonumber\\
&=&\alpha(T(Y,X),Z)+\beta(Y,X,Z)\,.\nonumber
\end{eqnarray}
This is the second identity in (\ref{cyc beta}). In view of this we also get
\begin{eqnarray}
&&\alpha(T(X,Y),Z)+\alpha(T(Y,Z),X)+\alpha(T(Z,X),Y)\nonumber\\
&=&\beta(Y,X,Z)-\beta(X,Y,Z)+\beta(Z,Y,X)\nonumber\\
&&-\beta(Y,Z,X)+\beta(X,Z,Y)-\beta(Z,X,Y)\nonumber\\
&=&0\nonumber
\end{eqnarray}
and this is (\ref{cyc beta2}). 
\end{proof}
\section{Lagrangian submanifolds in nearly K\"ahler six-manifolds}\label{sec 3}

By a well known theorem of Ejiri \cite{Ejiri-1981} Lagrangian submanifolds of $S^6$ are minimal.
In this section we will see that this is a special case of a much more general theorem.
As a consequence of Lemma \ref{lemma cyc} we will observe that Lagrangian
submanifolds in strict nearly K\"ahler manifolds of dimension six 
are minimal.
As the next theorem shows, the situation in three dimensions is quite special.
\begin{theorem}\label{Minimal_th_six}
Let $L^3$ be a Lagrangian immersion in a strict nearly K\"ahler six-manifold $M^6$. 
Then we have
\begin{eqnarray}
\alpha&=&0\label{3 alpha}\,,\\
\overrightarrow H&=&0\,.\label{3 mean}
\end{eqnarray}
In particular, any Lagrangian immersion in a strict nearly K\"ahler six-manifold is orientable and minimal.
\end{theorem}
\begin{proof}
Let $\{e_1,e_2,e_3\}$ be an orthonormal basis of $T_pL$ for a fixed point $p\in L$.
{F}rom the skew-symmetry of $\langle T(X,Y),Z\rangle$ we see that there exists a (nonzero) constant
$a$ such that $T(e_1,e_2)=ae_3$.
Then we also have $T(e_2,e_3)=ae_1,\, T(e_3,e_1)=ae_2$. The symmetry of $C$ implies
\begin{eqnarray}
\alpha(e_1,e_1)
&=&C(e_1,e_1,T(e_1,e_1))+C(e_2,e_1,T(e_2,e_1))+C(e_3,e_1,T(e_3,e_1))\nonumber\\
&=&0-aC(e_2,e_1,e_3)+aC(e_3,e_1,e_2)=0\nonumber
\end{eqnarray}
and
\begin{eqnarray}
\alpha(e_1,e_2)
&=&C(e_1,e_1,T(e_1,e_2))+C(e_2,e_1,T(e_2,e_2))+C(e_3,e_1,T(e_3,e_2))\nonumber\\
&=&aC(e_1,e_1,e_3)+0-aC(e_3,e_1,e_1)=0\,.\nonumber
\end{eqnarray}
Similarly we prove that $\alpha(e_i,e_j)=0$ for all $i,j=1,\dots,3$. This shows $\alpha=0$. But then (\ref{cyc mean})
also implies $\overrightarrow H=0$. 
The observation, that the frame $\{ e_1,e_2,e_3\}$ defines an orientation on
$L,$ finishes the proof.
\begin{remark}\label{bem type}
The constant $a$ in the formula $T(e_1,e_2)=ae_3$ from above is related to the type
constant $\alpha$ of the nearly K\"ahler manifold $M$ by the formula
$$a^2=\alpha\,.$$
This follows directly from the defining equation for the type constant which is
$$\langle\nabla_X(J)Y,\nabla_X(J)Y\rangle
=\alpha\left(|X|^2\cdot|Y|^2-\langle X,Y\rangle^2-\langle X,JY\rangle^2\right)\,.$$
In addition, recall that the type constant is related to the Ricci curvature (of the Levi-Civita connection)
of $M$ by
$$\operatorname{Ric}=5\alpha g\,,$$
so that the scalar curvature $s$ is $s=30\alpha$. In particular, strict six-dimensional
nearly K\"ahler manifolds are Einstein with $s>0$.
\end{remark}
\end{proof}

\section{The splitting theorem}\label{sec 5}
The following example shows that theorem \ref{Minimal_th_six} does not extend to eight dimensions:
\begin{example}
Let $L'\subset M^6_{SNK}$ be a (minimal) Lagrangian submanifold in a strict nearly K\"ahler
manifold $M^6_{SNK}$ and suppose $\gamma\subset\Sigma$ is a curve on a Riemann surface
$\Sigma$. Then the Lagrangian submanifold $L:=L'\times\gamma\subset M$ in
the nearly K\"ahler manifold $M:=\Sigma\times M_{SNK} $ is minimal, if and only if $\gamma$ is
a geodesic in $\Sigma$.
\end{example}
In this section we will see that this is basically the only counterexample to theorem \ref{Minimal_th_six}
that occurs in dimension eight.

Nearly K\"ahler manifolds $(M,J,g)$ split locally into a K\"ahler factor and a
strict nearly K\"ahler factor and under the assumption, that  $M$ is complete 
and simply connected this splitting is global \cite{Nagy-2002}. The natural
question answered in the following theorem is in which way Lagrangian
submanifolds lie in this decomposition. 
 
\begin{theorem}\label{split_th}
Let $(M,J,g)$ be a nearly K\"ahler manifold and $L\subset M$ a Lagrangian
submanifold. Then $M$ and $L$ split (locally) into $M=M_K\times M_{SNK}$, $L=L_K\times L_{SNK}$,
where $M_K$ is K\"ahler, $M_{SNK}$ is strict nearly K\"ahler and $L_K\subset M_K$,
$L_{SNK}\subset M_{SNK}$ are both Lagrangian. Moreover, the dimension of $L_K$ is given by
$$\operatorname{dim}(L_K)=\frac{1}{2}\operatorname{dim}(\operatorname{ker}(r))\,.$$
If the splitting of $M$ is global and $L$ is simply connected, then the splitting of $L$ 
is global as well.
\end{theorem}
\begin{proof}~

\begin{enumerate}[i)]
\item\label{p1}
We define
$$K_p:=\{X\in T_pM:rX=0\}\,,\quad K_p^\perp:=\{Y\in T_pM:\langle X,Y\rangle=0\,,\forall\, X\in K_p\}\,.$$
Because of $\nb r=0, \nb g=0$ this defines two orthogonal smooth distributions
$$\mathscr{D}_K:=\bigcup_{p\in M}K_p\,,\quad \mathscr{D}_{SNK}:=\bigcup_{p\in M}K_p^\perp$$
on $M$. 
\item\label{p2}
Since 
$$\langle rX,Y\rangle=\sum_{i=1}^{2n}\langle\nabla_X(J)e_i,\nabla_Y(J)e_i\rangle
=\langle\nabla_X(J),\nabla_Y(J)\rangle_{\operatorname{End}(TM)}$$
$r$ is a selfadjoint and positive semidefinite operator. Then
$\langle rX,X\rangle=||\nabla_X(J)||^2$ implies 
$$X\in \operatorname{ker}(r)\quad\Leftrightarrow\quad \nabla_X(J)=0\,.$$
\item\label{p3}
Let $K\in\Gamma(\mathscr{D}_K)$, $X\in\mathscr{X}(M)$ be smooth vector fields. Then $\nb r=0$ and $r(K)=0$ imply
$$X(r(K))=0=\nb_X(r)K+r(\nb_XK)=r(\nb_XK)\,.$$
Hence $\nb_XK\in\Gamma(\mathscr{D}_K)$ as well.
\item
Suppose $K_1,K_2\in\Gamma(\mathscr{D}_K)$. Then
\begin{eqnarray}
[K_1,K_2]
&=&\nabla_{K_1}K_2-\nabla_{K_2}K_1\nonumber\\
&=&\nb_{K_1}K_2+\frac{1}{2}J\nabla_{K_1}(J)K_2-\nb_{K_2}K_1-\frac{1}{2}J\nabla_{K_2}(J)K_1\nonumber\\
&\overset{(\ref{p2})}{=}&\nb_{K_1}K_2-\nb_{K_2}K_1\overset{(\ref{p3})}{\in}\Gamma(\mathscr{D}_K)\,.\nonumber
\end{eqnarray}
Hence $\mathscr{D}_K$ is an integrable distribution.
\item
The orthogonal distribution $\mathscr{D}_{SNK}$ is also integrable. For $X\in\mathscr{X}(M)$, $S\in\Gamma(\mathscr{D}_{SNK})$, $K\in\Gamma(\mathscr{D}_K)$ we compute
\begin{eqnarray}
X\langle K,S\rangle&=&0\nonumber\\
&=&\langle\nabla_XK,S\rangle+\langle K,\nabla_XS\rangle\nonumber\\
&=&\langle\nb_XK+\frac{1}{2}J\nabla_X(J)K,S\rangle+\langle K,\nabla_XS\rangle\nonumber\\
&\overset{(\ref{p3})}{=}&\frac{1}{2}\langle J\nabla_X(J)K,S\rangle+\langle K,\nabla_XS\rangle\nonumber\\
&\overset{(\ref{nk4})}{=}&-\frac{1}{2}\langle J\nabla_K(J)X,S\rangle+\langle K,\nabla_XS\rangle\nonumber\\
&\overset{(\ref{p2})}{=}&\langle K,\nabla_XS\rangle\,,\nonumber
\end{eqnarray}
so that $\nabla_XS\in\Gamma(\mathscr{D}_{SNK})$. Then again
$$[S_1,S_2]=\nabla_{S_1}S_2-\nabla_{S_2}S_1\in\Gamma(\mathscr{D}_{SNK})\,.$$
\item
The splitting theorem of de Rham can now be applied to the distributions $\mathscr{D}_K$ and
$\mathscr{D}_{SNK}$ and the nearly K\"ahler manifold $(M,J,g)$ splits
(locally) into a Riemannian product
$$(M,J,g)=(M_K,J_K,g_K)\times (M_{SNK},J_{SNK},g_{SNK})\,,$$
where $TM_K=\mathscr{D}_K$, $TM_{SNK}=\mathscr{D}_{SNK}$ and where $(M_K,J_K,g_K)$ is
K\"ahler and $(M_{SNK},J_{SNK},g_{SNK})$ is strict nearly K\"ahler.
\item
Now let $L\subset M=M_K\times M_{SNK}$ be Lagrangian. We prove that $r$ leaves
tangent and normal spaces of $L$ invariant. To see this, we fix an adapted local orthonormal
frame field $\{e_1,\dots,e_{2n}\}$ of $M$ such that $e_1,\dots,e_n$ are tangent
to $L$ and $e_{n+1}=Je_1,\dots,e_{2n}=Je_n$ are normal to $L$. Since for any three 
vectors $X,Y,Z$ we have
$$\langle\nabla_X(J)JZ,\nabla_Y(J)JZ\rangle=\langle J\nabla_X(J)Z,J\nabla_Y(J)Z\rangle
=\langle\nabla_X(J)Z,\nabla_Y(J)Z\rangle\,,$$
we obtain
\begin{eqnarray}
\langle rX,Y\rangle
&=&\sum_{i=1}^{2n}\langle\nabla_X(J)e_i,\nabla_Y(J)e_i\rangle\nonumber\\
&=&\sum_{i=1}^{n}\langle\nabla_X(J)e_i,\nabla_Y(J)e_i\rangle
+\sum_{i=1}^{n}\langle\nabla_X(J)Je_i,\nabla_Y(J)Je_i\rangle\nonumber\\
&=&2\sum_{i=1}^{n}\langle\nabla_X(J)e_i,\nabla_Y(J)e_i\rangle.\nonumber
\end{eqnarray}
Now, if $X\in TL$, $Y\in T^\perp L$, then by Lemma \ref{lemma 1} we have
$$\nabla_X(J)e_i\in T^\perp L\,,\quad\nabla_Y(J)e_i\in TL\,,$$
so that
$$\langle\nabla_X(J)e_i,\nabla_Y(J)e_i\rangle=0\,,\quad\forall\, i=1,\dots,n\,.$$
{F}urther it follows 
$$\langle rX,Y\rangle=2\sum_{i=1}^n\langle\nabla_X(J)e_i,\nabla_Y(J)e_i\rangle=0\,.$$
Since this works for any $X\in TL, Y\in T^\perp L$ and since $r$ is selfadjoint
we conclude
$$r(TL)\subset TL\,,\quad r(T^\perp L)\subset T^\perp L\,.$$
At a given point $p\in L$ we may now choose an orthonormal basis $\{f_1,\dots,f_n\}$
of $T_pL$ that consists of eigenvectors of $r_{|TL}$ considered as an endomorphism of $TL$.
Since $[r,J]=0$ and $L$ is Lagrangian, the set $\{f_1,\dots,f_n,Jf_1,\dots,Jf_n\}$ then 
also determines an orthonormal eigenbasis of $r\in\operatorname{End}(TM)$. In particular,
since $J$ leaves the eigenspaces invariant, $K_p=\operatorname{ker}(r(p))$ and $T_pL$
intersect in a subspace $K^L_p$ of dimension $\frac{1}{2}\dim(K_p)=\frac{1}{2}\dim(M_K)$.
{F}or the same reason $K_p^\perp\cap T_pL$ gives an $\frac{1}{2}\dim(M_{SNK})$-dimensional
subspace. The corresponding distributions, denoted by $\mathscr{D}_K^L$ and 
$\mathscr{D}_{SNK}^L$ are orthogonal and both integrable, since in view of
$$\mathscr{D}_K^L=\mathscr{D}_K\cap TL\,,\quad\mathscr{D}_{SNK}^L=\mathscr{D}_{SNK}\cap TL$$
they are given by intersections of integrable distributions. We may now apply again the
splitting theorem of de Rham to the Lagrangian submanifold. This completes the proof.
\end{enumerate}
\end{proof}

A more detailed analysis of the proof of the last theorem shows

\begin{corollary}\label{cor r}
If  $L\subset M$
is Lagrangian and $p\in L$ a fixed point, then to each eigenvalue $\lambda$ of
the operator $r$ at $p$ there exists a basis $e_1,\dots,e_k,f_1,\dots,f_k$ of eigenvectors
of $\operatorname{Eig}(\lambda)$ 
such that $e_1,\dots,e_k\in T_pL$, $f_1,\dots,f_k\in T_p^\perp L$.
Here, $2k$ denotes the multiplicity of $\lambda$.
\end{corollary}

Let us recall the situation in dimension eight and ten \cite{Gray-1976,Nagy-2002}.
\begin{proposition} \label{NK_dim8_10}~

\begin{enumerate}[(i)]
\item
Let $M^8$ be a simply connected complete nearly K\"ahler manifold of dimension eight. 
  Then $M^8$ is a Riemannian
  product $M^8= \Sigma \times M_{SNK}^6$ of a Riemannian surface $\Sigma$ and  a
  six-dimensional strict nearly K\"ahler manifold $M_{SNK}^6.$
\item
Let $M^{10 }$ be a simply connected complete nearly K\"ahler manifold of
dimension ten. Then $M^{10}$ is either the product $M_K^4 \times M_{SNK}^6$ of a K\"ahler surface $M_K^4$ and a
  six-dimensional strict nearly K\"ahler manifold $M_{SNK}^6$ or $M$ is a
  twistor space over a positive, eight dimensional quaternionic K\"ahler manifold.
\end{enumerate}
\end{proposition}
Note, that any complete, simply connected  eight dimensional quaternionic manifold
equals one of the following three spaces: $\mathbb{HP}^2, \mathbb{G}r_2(\complex{2}), G_2/SO(4)$

In the next theorem, part (i) and (ii) collect
the information on Lagrangian submanifolds in nearly K\"ahler manifolds of
dimension eight and ten. 

\begin{theorem}\label{theoC}~

\begin{enumerate}[(i)]
\item
Let $L$  be a Lagrangian submanifold in a simply connected nearly K\"ahler manifold $M^8$.
Then $M^8=\Sigma \times M^6_{SNK}$, where
$\Sigma$ is a Riemann surface, $M^6_{SNK}$ is strict nearly K\"ahler and $L=\gamma\times L'$ is a
  product of a (real) curve $\gamma\subset\Sigma$ and a  minimal Lagrangian submanifold $L' \subset  M^6_{SNK}.$
\item\label{theoC2}
Let $L$  be a Lagrangian submanifold in a simply
  connected complete nearly K\"ahler manifold $M^{10}$, then either
  \begin{enumerate}[(a)]
  \item $M^{10}=M_K^4 \times M^6_{SNK}$ and the manifold  $L= S \times L'$ is a product of a Lagrangian (real) surface $S\subset
    M_K^4$ and a minimal Lagrangian submanifold $L' \subset  M^6_{SNK}$ or
  \item\label{theoC2b} the manifold $L$ is a  Lagrangian submanifold in a  twistor space
    over a positive, eight dimensional quaternionic K\"ahler manifold.
 \end{enumerate}
\item 
Let $M_1$, $M_2$ be two nearly K\"ahler manifolds. Denote the operator $r$ on $M_i$ by $r_i$,
$i=1,2$. 
If $\operatorname{Spec}(r_1)\cap\operatorname{Spec}(r_2)=\emptyset$ and $L\subset M_1\times M_2$
is Lagrangian, then $L$ splits (locally) into $L=L_1\times L_2$, where $L_i\subset M_i$, $i=1,2$
are Lagrangian. If $L$ is simply connected, then the decomposition is global.
\end{enumerate}
\end{theorem}
\begin{proof}
This is a combination of the results of Theorem \ref{split_th}, Corollary \ref{cor r}
and Proposition \ref{NK_dim8_10}.
\end{proof}
Theorem \ref{theoC} (\ref{theoC2}), part (b) motivates the discussion of Lagrangian
submanifolds in twistor spaces in the subsequent section. Indeed, the results derived in 
the next section imply that Lagrangian submanifolds in twistor spaces are, regardless
their dimension, always minimal.
\section{Lagrangian submanifolds in twistor spaces}\label{twistor_sec}

As mentioned in the introduction an important class of examples for nearly K\"ahler 
manifolds is given by
 twistor spaces $Z^{4n+2}$ over positive quaternionic K\"ahler manifolds
 $N^{4n}.$ We recall that the twistor space is the bundle of almost complex
 structures over $N.$  It can be endowed with a K\"ahler  structure $(Z,J^Z,g^Z),$ such that the
 projection $\pi \,:\, Z \rightarrow N$ is a  Riemannian submersion with totally
 geodesic fibers $S^2.$ Denote by $\mathcal H$ and $\mathcal V$ the horizontal
 and the vertical distributions of the submersion $\pi.$ Then the direct
 decomposition 
\begin{equation}
TZ = \mathcal H \oplus \mathcal V \label{TM_dec_tw}
\end{equation}
 is orthogonal and compatible
 with the complex structure $J^Z.$ Let us consider now a second almost
 hermitian structure $(J,g)$ on $Z$ which is defined by
$$ g:= 
\begin{cases} 
 g^Z(X,Y), \mbox{ for } X,Y \in \mathcal{H}, \\
 \frac{1}{2} g^Z(V,W), \mbox{ for } V,W \in \mathcal{V},\\
g^Z(V,X)=0,   \mbox{ for } V \in \mathcal{V},
X \in \mathcal{H} 
\end{cases}$$
and
$$ J:= 
\begin{cases} 
 J^Z \mbox{ on } \mathcal{H}, \\ 
-J^Z \mbox{ on } \mathcal{V}.
\end{cases}$$

Note, that in view of (\ref{TM_dec_tw}), the decomposition $TZ=\mathcal H\oplus\mathcal V$ is also compatible w.r.t. $J$ and orthogonal w.r.t. $g$.
Then the following result is well-known.
\begin{proposition} \label{can_NK_tw}
The manifold $(Z,J,g)$ is a strict nearly K\"ahler manifold 
and the distributions $\mathcal V$ and $\mathcal H$ are 
parallel w.r.t. the connection $\nb.$ 
The projection $\pi$ is also a Riemannian submersion 
with totally geodesic fibers for the metric $g.$
\end{proposition}
We summarize some information which will be useful in the later text.
\begin{lemma} \label{info_lemma}
In the situation of Proposition \ref{can_NK_tw} we have
\begin{enumerate}[(a)]
\item
The torsion $T =-
J\nabla (J)$ of the canonical connection satisfies
\begin{eqnarray}
T(X,Y) \in \mathcal V,&& \mbox{ for } X,Y \in \mathcal H, \label{Tor_1_equ}\\
T(X,V) \in \mathcal H,&& \mbox{ for } X \in \mathcal H, V \in \mathcal V,\label{Tor_2_equ}\\
T(U,V) =0, && \mbox{ for } U,V \in \mathcal V. \label{Tor_3_equ}
\end{eqnarray}
\item
The association 
\begin{equation}\label{ass 1}
\mathcal H \ni X \mapsto T(Y,X)\in\mathcal V
\end{equation}
is surjective for $0\ne Y \in \mathcal H$ and the map
\begin{equation}\label{ass 2}
\Phi^V\,:\, \mathcal H \ni X \mapsto T(V,X)
\end{equation}
with $0\ne V \in \mathcal V$ is invertible and squares
to $-\kappa^2 Id_{\mathcal H}$ for some $\kappa \in \mathbb{R}.$
\item
The operator $r$ has eigenvalues $\lambda_{\mathcal{H}}=4\kappa^2$, $\lambda_{\mathcal V}=4n\kappa^2$. 
If $n>1$, then the eigenbundle of $\lambda_{\mathcal H}$ is $\mathcal{H}$ and $\mathcal{V}$ 
is the eigenbundle of $\lambda_{\mathcal V}$.
\end{enumerate}
\end{lemma}
\begin{proof} Property \eqref{Tor_1_equ} follows from Propostion 4.1 (ii) of
  \cite{Nagy-2002}. The second property \eqref{Tor_2_equ} follows, since
  $g(T(\cdot,\cdot),\cdot)$ is a three-form and \eqref{Tor_3_equ} follows by a direct computation from the nearly
K\"ahler condition. The second part of (b) is a consequence of Proposition 4.1
(i) of \cite{Nagy-2002}. 
We claim that for a fixed $Y \in
\mathcal H$ there exists $X \in \mathcal H,$ such that $\mathcal V \ni V=T(Y,X) \ne 0.$
This yields $JV=-T(Y,JX) \ne 0$ and proves the first part of (b). It remains to
prove the claim:  $T(Y,X)=0$ for all $X \in \mathcal H$ implies $g(T(Y,X),V)=-
g(J\nabla_Y(J)X,V)=g(\nabla_Y(J)X,JV)=-g(\nabla_Y(J)V,JX) =0$ for all $X \in \mathcal H$
and $V \in \mathcal V.$ In other words $\nabla_Y(J) =0$ in contradiction to strict
nearly K\"ahler. Part (c) can be found in Corollary 4.1 (i) of \cite{Nagy-2002}.
\end{proof}
In the rest of this section we consider a nearly K\"ahler manifold $(M=Z,J,g)$
of this type and study Lagrangian submanifolds $L\subset M.$ 

\begin{lemma} \label{2nd_fund_tw}
Let $L^{2n+1}\subset M^{4n+2}$ be a Lagrangian submanifold in a twistor space as 
in Proposition \ref{can_NK_tw}. 
Then the second
fundamental form $II$ satisfies
\begin{eqnarray}
&II(X,Y) \in  \pi^{\mathcal H}(T^\perp L),& \mbox{ for } X,Y \in  \pi^{\mathcal H}(TL), \\
&II(X,Y) \in  \pi^{\mathcal V}(T^\perp L),& \mbox{ for } X,Y \in  \pi^{\mathcal V}(TL),\\
&II(X,Y) =0,& \mbox{ for } X \in  \pi^{\mathcal H}(TL), Y \in  \pi^{\mathcal V}(TL).
\end{eqnarray}
Here $\pi^{\mathcal H}(TL)$ and $\pi^{\mathcal V}(TL)$ are the orthogonal projections of $TL$ to
$\mathcal H$ resp. $\mathcal V$.
\end{lemma}
\begin{proof}
The second fundamental form is given by 
$C(X,Y,Z)=\langle \nb_X Y, JZ \rangle$ 
for $X,Y,Z \in TL.$ The lemma follows since
the decomposition \eqref{TM_dec_tw} is $\nb$-parallel, orthogonal and
$J$-invariant and as the tensor $C$ is completely symmetric.
\end{proof}
\begin{remark}
As will be shown in the next theorem, for $n>1$ we have
$$\pi^{\mathcal H}(TL)=\mathcal H\cap TL\,,\quad\pi^{\mathcal V}(TL)=\mathcal V\cap TL\,.$$
\end{remark}

\begin{theorem} \label{min__Lgra_tw}
Let $L^{2n+1}\subset M^{4n+2}$ be a Lagrangian submanifold in a twistor space as 
in Proposition \ref{can_NK_tw}. Then
$L$ is minimal. If $n>1$, then the tangent space of $L$ splits into a one-dimensional
vertical part and a $2n$-dimensional horizontal part. Moreover, the
second fundamental form $II$ of the vertical normal direction vanishes completely
if $n>1$. 
\end{theorem}
\begin{proof}~

\begin{enumerate}[i)]
\item
By theorem \ref{Minimal_th_six} it suffices to consider the case $n>1$.
\item\label{Lemma_distr_D}
Let $L\subset M$ be a Lagrangian submanifold. Since $n>1$, the two eigenvalues
$\lambda_{\mathcal H},\,\lambda_{\mathcal V}$ of $r$ are distinct and the eigenspace $\mathcal{V}$ of
$\lambda_{\mathcal V}$ is two-dimensional. By Corollary \ref{cor r} this induces a
one-dimensional vertical tangential distribution $\mathcal D$ on $L$.
Then, by the Lagrangian condition, $\mathcal D := \pi^{\mathcal V}(TL) $.
\item
Denote by 
$\mathcal D^\perp$ the orthogonal complement of $\mathcal D$ in $TL.$
The trace of the second fundamental form $II$ of $L$ restricted to
$\mathcal D^\perp$ is zero.
\begin{proof}
{F}irst we observe that by Lemma \ref{2nd_fund_tw}  we can restrict the
second fundamental form $II$ to $\mathcal D^\perp =\pi^{\mathcal H}(TL).$ 
We fix
$U \in \mathcal D$ of unit length.  Using Lemma \ref{lemma 1} and Lemma \ref{info_lemma} we observe, 
that $\Phi(X):= \frac{1}{\kappa} J\nabla_U(J)X$ defines an (almost) complex
structure on $\mathcal D^\perp$ which is compatible with the metric. 
{F}irst we observe with $X \in \mathcal D^\perp.$
\begin{eqnarray*}
II(X,X)&=&-II(X,\Phi(\Phi(X)))=-\frac{1}{\kappa}II(X,J\nabla_U(J)\,\Phi(X))\\
&=&-\frac{1}{\kappa}J\left[\nabla_{II(X,U)}(J)\Phi(X)
+\nabla_U(J)II(X,\Phi(X)) \right]\\
&=&-\frac{1}{\kappa}J\nabla_U(J)\, II(X,\Phi(X))=-\Phi II(X,\Phi(X))\,. 
\end{eqnarray*}
After polarizing we obtain
\begin{equation}\label{pluriminimal}
II(\Phi X,\Phi Y)=-II(X,Y)\,,\quad\forall\, X,\,Y\in\mathcal D^\perp\,.
\end{equation}
In particular, taking a trace over (\ref{pluriminimal}) we get
\begin{eqnarray*}
\operatorname{tr}^{\mathcal D^\perp}II &=& 0\,.
\end{eqnarray*}
\end{proof}
\item
We have
$$\alpha(X,Y)=0\,,\quad\forall\, X,\,Y\in\mathcal D^\perp\,.$$
\begin{proof}
By (ii) we may choose an orthonormal frame $\{e_1,\dots,e_{2n+1}\}$
of $TL$ such that $e_1,\dots,e_{2n}\in\mathcal D^\perp$ and $e_{2n+1}\in\mathcal D$.
Since
$$\alpha(X,Y)=\sum_{i=1}^{2n+1}C(e_i,X,T(e_i,Y))$$
and the tensor $C$ is fully symmetric we see that by Lemma \ref{2nd_fund_tw} all terms on the RHS vanish since either $e_i\in \mathcal D=\pi^{\mathcal V}(TL)$, $X\in\mathcal D^\perp=\pi^{\mathcal H}(TL)$ or $e_i, X\in\mathcal D^\perp$ and  $T(e_i,Y)\in\mathcal D$ (cf. Lemma \ref{info_lemma}).
\end{proof}
\item
By (iii) the mean curvature vector $\overrightarrow H$ satisfies
$$\overrightarrow H\in\mathcal D\,.$$
{F}rom (\ref{cyc mean}) and (iv) we get
\begin{equation}\label{mean 1}
\langle J\overrightarrow H, T(X,Y)\rangle=0\,,\quad\forall\,X,\,Y\in\mathcal D^\perp\,.
\end{equation}
Since $J$ maps $\mathcal V$ to itself, we also have $J\overrightarrow H\in\mathcal D\subset\mathcal V$.
Now we choose $X\in \mathcal{D}^\perp$ and $\tilde Y\in\mathcal H$ with
$$T(X,\tilde Y)=J\overrightarrow H\,.$$
This is possible since the map $\mathcal H\ni \tilde Y\mapsto T(X,\tilde Y)\in\mathcal V$ is 
surjective by (\ref{ass 1}).
Let $\tilde Y=Y+Y^\perp$ be the orthogonal decomposition of $\tilde Y$ into the tangent
and normal parts of $\tilde Y$. Note that $Y,Y^\perp$ are both horizontal. We have
$$T(X,Y^\perp)=0$$
since $T(X,\cdot)$ maps tangent to tangent and normal to normal vectors 
and
$$T(X,\tilde Y)=T(X,Y)+T(X,Y^\perp)=J\overrightarrow H\in TL\,.$$
Therefore there exist two tangent vectors $X,Y\in\mathcal{D}^\perp$ with 
$$T(X,Y)=J\overrightarrow H\,.$$
This implies
$$|\overrightarrow H|^2=\langle J\overrightarrow H,T(X,Y)\rangle\overset{(\ref{mean 1})}{=}0$$
which proves that the mean curvature vector vanishes.
{F}rom this, the fact that $\mathcal D$ is one-dimensional and from Lemma \ref{2nd_fund_tw} it follows
that $II(V,\cdot)=0$ for any $V\in\mathcal D$.
\end{enumerate}
\end{proof}

\begin{corollary}
Let $L\subset M$ be a Lagrangian submanifold in a twistor space $M^{4n+2}$ as above with $n>1$. 
Then the integral manifolds $c$ of the distribution $\mathcal D$  are geodesics
(hence locally great circles) in the totally geodesic fibers  $S^2.$ 
 \end{corollary}
\begin{proof}
The last theorem implies that the geodesic curvature 
vanishes and that in consequence an integral manifold $c$ of $\mathcal D$ is totally geodesic in the fibers.
\end{proof}

\section{Invariant Lagrangian submanifolds of three-symmetric spaces}\label{sec 7}
The idea of a three-symmetric space is to replace the symmetry of order two as
in the case of a symmetric space by a symmetry  of order three. 
Nearly K\"ahler geometry on such spaces was first studied  in \cite{Gray-1972,Gray-Wolf}. \\
Like symmetric spaces three-symmetric spaces have a homogenous model, which
we shortly resume. Let $G$ be a connected Lie group and $s$ an automorphism of
order $3$ and let $G^s_0 \subset H \subset G^s$ be a subgroup contained in
the fix-point set $G^s$ of $s.$ The differential $s_*$ decomposes 
$$ \mathfrak g \otimes \mathbb C = \mathfrak h \otimes \mathbb C \oplus
\mathfrak m^+  \oplus \mathfrak m^-$$
into the eigenspaces of $s_*$ with eigenvalues $1$ and $\frac{1}{2}(-1 \pm
\sqrt{-3}).$ With the definition $\mathfrak m := (\mathfrak m^+  \oplus
\mathfrak m^-)\cap \mathfrak g$ the decomposition 
\begin{equation} \label{3-sym_dec} \mathfrak g=\mathfrak h \oplus \mathfrak  m  \end{equation}
is reductive and $G/H$ is a reductive homogenous space. The canonical complex
structure is then defined by
$$ s_*| \mathfrak m= -\frac{1}{2} Id + \frac{\sqrt{3}}{2} J.$$
The choice of an $Ad(H)$-invariant and $s_*$-invariant (pseudo-)metric $B$ on $\mathfrak m$
makes $G/H$ into a (pseudo-)Riemannian three-symmetric space, such that $B$ is
almost hermitian with respect to $J.$ \\
By Proposition 5.6 of \cite{Gray-1972} the data $(M=G/H,J,B)$ defines a nearly K\"ahler manifold if and only if the
decomposition \eqref{3-sym_dec} is naturally reductive, i.e. $$ B([X,Y]_{\mathfrak m},Z)=B(X,[Y,Z]_{\mathfrak m})\mbox{    for } X,Y,Z \in \mathfrak m.$$
Since $M$ is homogenous and naturally reductive the Levi-Civita connection in
the canonical base point $p$ (after identifying $\mathfrak m$ and the tangent space) is given by
\begin{equation} \label{LC_3_sym_NK} 2B(\nabla_XY,Z) =  B([X,Y],Z) \mbox{
    for } X,Y,Z \in \mathfrak m.\end{equation}
Moreover it is $B(\nabla_X(J)Y,JZ)= B(X,[Y,Z]) \mbox{
    for } X,Y,Z \in \mathfrak m$ which implies $$J \nabla_X(J)Y =
[X,Y]_{\mathfrak m}\mbox{    for } X,Y \in \mathfrak m $$ in $p.$ With this
information the canonical connection $\nb$ and the Levi-Civita connection $\nabla$  are related in $p$ by \begin{equation} \label{nb_3_sym_NK}
\nb_XY - \nabla_XY = -\frac{1}{2} J \nabla_X(J)Y= -\frac{1}{2}[X,Y]_{\mathfrak m}\mbox{
    for } X,Y,Z \in \mathfrak m. \end{equation}

\begin{theorem}
Let $L$ be an invariant  Lagrangian submanifold in $M=G/H$ with the above given nearly
K\"ahler structure. Then $L$ is totally geodesic. In particular $L$ is
minimal. 
\end{theorem}
\begin{proof}
By Proposition \ref{2nd_fund_info} the second fundamental form is given by  
$$B(II(X,Y),JZ)= B(\nabla_XY,JZ) = B(\nb_XY,JZ)\,.$$ 
This shows with equation \eqref{nb_3_sym_NK}
$$0=B(\nb_XY-\nabla_XY,JZ)= -\frac{1}{2}B([X,Y]_{\mathfrak m},J Z) \mbox{
    for } X,Y,Z \in \mathfrak m\cap T_pL.$$
Using equation \eqref{LC_3_sym_NK} we obtain 
$$B(II(X,Y),JZ) = B(\nabla_XY,JZ)= \frac{1}{2}  B( [X,Y]_{\mathfrak
  m},JZ),$$
which is zero by the previous equation. This shows that $II=0,$ i.e. $L$ is totally geodesic and consequently minimal.
\end{proof}
\subsection{Examples on $S^3\times S^3$}
A convenient way to introduce the nearly K\"ahler structure on 
$M=S^3\times S^3$ is to identify $S^3$ with $SU(2)$ and $M$ with $SU(2) \times
SU(2).$ Let us denote by $\mathcal E =(e_1,e_2,e_3)$ a frame of left invariant one-forms
on the first factor and  by $\mathcal F= (f_1,f_2,f_3)$ a frame left invariant one-forms
on the second factor. One can choose the frames $\mathcal E$ and $\mathcal F$
(see for example \cite{Butruille-2005}) in such a way that the fundamental two-form $\omega$ is given by
\begin{equation}
\omega = c\sum_{i=1}^3 e_i \wedge f_i \mbox{ with } c\in \mathbb{R}.
\end{equation}
With this preparation we characterize left invariant Lagrangian submanifolds
in the following theorem.

\begin{theorem}
Let $L\subset M=S^3\times S^3$ be an invariant Lagrangian submanifold. 
Then $L$ is $S^3$ after inclusion in the first or second factor or  
an integral manifold of the left invariant distributions 
$\mathcal D= \{ X \oplus DX \,|\,
X \in \mathfrak{su}(2)\} \subset\mathfrak{su}(2) \oplus \mathfrak{su}(2) $ with
$D\in \{\text{diag}(1,1,1), \text{diag}(1,1,-1), \text{diag}(1,-1,-1), \text{diag}(-1,-1,-1)\}.$
\end{theorem}
\begin{proof}
After the action of an element of $SU(2) \times SU(2)$ we suppose, that the
neutral element $e$ lies in $L$ and consider
$\mathcal L=T_eL \subset T_e(SU(2) \times SU(2)) = \mathfrak{su}(2) \oplus  \mathfrak{su}(2).$ It is not difficult 
to see that inclusion of $S^3$ in the first or second factor yields a
Lagrangian submanifold. Therefore we suppose that $T_eL$ is 
not one of the factors. Since $\mathcal L$ is Lagrangian we can express it as
a graph, i.e. $\mathcal L= X \oplus AX$ for $X \in \mathfrak{su}(2)$ and a linear
map $A\,:\, \mathfrak{su}(2) \rightarrow \mathfrak{su}(2).$ The integrability
of $\mathcal L$ yields that $[X \oplus AX,Y \oplus AY]= [X,Y] \oplus [AX,AY]$
equals $[X,Y] \oplus A[X,Y].$ This is equivalent to $A[X,Y]=[AX,AY].$ If we
identify  $\mathfrak{su}(2)$ with $\mathbb{R}^3$ endowed with the cross
product, then the non-trivial solutions of this constraint are $A^tA=1.$
The condition $\omega_{| \mathcal L}=0$ is equivalent to $A^t=A.$ This means
we end up with a symmetric $3\times 3$ matrix satisfying $A^2=1.$ Using
Sylvester's theorem the solutions are $SO(3)$-orbits of one of diagonal
matrices
$D\in \{ \text{diag}(1,1,1), \text{diag}(1,1,-1), \text{diag}(1,-1,-1), \text{diag}(-1,-1,-1) \}.$
The action of $SO(3)$ corresponds to the choice of some basis of
$ \mathfrak{su}(2)$ w.r.t. to some basis left invariant one-forms.
\end{proof} 

\section{Deformations of Lagrangian submanifolds in nearly K\"ahler manifolds}\label{sec 8}
Our aim in this section is to study the space of deformations of a given Lagrangian
(and hence minimal Lagrangian) submanifold $L$ in a strict six-dimensional
nearly K\"ahler manifold $M^6$. In an article by Moroianu, Nagy and Semmelmann
\cite{MNS-2008} the deformation space of nearly K\"ahler structures on six-dimensional
nearly K\"ahler manifolds has been related to the space of coclosed eigenforms
of the Hodge-Laplacian. As we will show below, a similar statement holds for the deformation
of Lagrangian submanifolds in strict nearly K\"ahler six-manifolds.

To this end we assume that
$$F:L\times(-\epsilon,\epsilon)\to M$$
is a smooth variation of Lagrangian immersions $F_t:=F(\cdot,t):L\to M$, $t\in(-\epsilon,\epsilon)$
into a nearly K\"ahler manifold $M$. Let 
$$V:=\frac{d}{dt} F_t$$
denote the variation vector field. Since tangential deformations correspond to
diffeomorphisms acting on $L$, we may assume w.l.o.g. that $V\in\Gamma(T^\perp L)$
is a normal vector field. The Cartan formula and $F^*_t\omega=0$ for all $t$
then implies 
$$0=d(i_V\omega)+i_Vd\omega$$
holds everywhere on $L$. By the nearly K\"ahler condition this is equivalent to
\begin{equation}\label{def 1}
d(V\cont\omega)+3V\cont\nabla\omega=0
\end{equation}
on $L$. Let us define the variation $1$-form $\theta\in\Omega^1(L)$ by
$$\theta:=V\cont\omega\,.$$
\begin{theorem}
Let $F_t:L\to M$ be a variation of Lagrangian immersions in a six-dimensional
nearly K\"ahler manifold $M$. Then the variation $1$-form $\theta$ is a coclosed
eigenform of the Hodge-Laplacian, where the eigenvalue $\lambda$ satisfies
$\lambda=9\alpha$ with the type constant $\alpha$ of $M$.
\end{theorem}
\begin{proof}
{F}or $X,Y\in TL$ and $V\in T^\perp L$ we compute
\begin{eqnarray}
(V\cont\nabla\omega)(X,Y)
&=&\nabla\omega(V,X,Y)\nonumber\\
&=&\langle \nabla_X(J)Y,V\rangle\nonumber\\
&=&\langle J\nabla_X(J)Y,JV\rangle\nonumber\\
&=&-\langle JV,T(X,Y)\rangle\,.\nonumber
\end{eqnarray}
Since $T$ induces an orientation on the Lagrangian submanifold by the three-form
$$\tau(X,Y,Z):=\langle T(X,Y),Z\rangle\,,$$
we obtain a naturally defined $*$-operator $*:\Omega^p(L)\to\Omega^{3-p}(L)$ which
for $1$-forms is given by
$$*\phi:=\frac{1}{\sqrt{\alpha}}\,\phi\circ T\,.$$
Here, $\alpha$ is the type constant of $M$ (cf. Remark \ref{bem type}).
This implies that equation (\ref{def 1}) can be rewritten in the form
\begin{equation}\label{def 2}
d\theta=3\sqrt{\alpha}\,{*\theta}\,.
\end{equation}
Consequently
$$\delta\theta=*d{*\theta}=0$$
and
\begin{eqnarray}
\delta d\theta
&=&3\sqrt{\alpha}\,{*d{**}}\theta\nonumber\\
&\overset{(*^2=\operatorname{Id})}{=}&3\sqrt{\alpha}\,{*d}\theta\nonumber\\
&\overset{(\ref{def 2})}{=}&9\alpha**\theta\nonumber\\
&=&9\alpha\theta\,.\nonumber
\end{eqnarray}
In total
$$\Delta_{\operatorname{Hodge}}\theta=(\delta d+d\delta)\theta=9\alpha\theta\,.$$
This proves the theorem. By remark \ref{bem type} this is equivalent to
$$\Delta_{\operatorname{Hodge}}\theta=\frac{3}{10}s\,\theta\,,$$
where $s$ is the scalar curvature of $M$.
\end{proof}

\bibliographystyle{plain}   
   
\end{document}